\documentclass[11pt]{article}
\usepackage{amsfonts,amscd,amssymb}
\newtheorem{definition}{Definition}[section]

\newtheorem{proposition}[definition]{Proposition}

\newtheorem{remark}[definition]{Remark}

\newtheorem{theorem}[definition]{Theorem}
\newtheorem{example}[definition]{Example}

\def\rawo\lonra{\longrightarrow}

\def\ot{\otimes}

\newcommand{\eqref}[1]{(\ref{eq:#1})}

\newenvironment{proof}{{\it Proof.}}{\hfill $ \square $ \vskip 4mm}
\begin{document}
\title{On some classes of lazy cocycles and categorical structures
\thanks{Research partially supported by the 
EC programme LIEGRITS,  
RTN 2003, 505078, and by the project ``New techniques in  
Hopf algebras and graded ring theory'' of the Flemish and Romanian 
Ministries of Research.}}
\author
{Florin Panaite\\
Institute of Mathematics of the 
Romanian Academy\\ 
PO-Box 1-764, RO-014700 Bucharest, Romania\\
e-mail: Florin.Panaite@imar.ro
\and 
Mihai D. Staic\thanks{Permanent address: Institute of Mathematics of the  
Romanian Academy, 
PO-Box 1-764, RO-014700 Bucharest, Romania.}\\
SUNY at Buffalo\\
Amherst, NY 14260-2900, USA\\
e-mail: mdstaic@buffalo.edu
\and
Freddy Van Oystaeyen\\
Department of Mathematics and Computer Science\\
University of Antwerp, Middelheimlaan 1\\
B-2020 Antwerp, Belgium\\
e-mail: Francine.Schoeters@ua.ac.be}
\date{}
\maketitle
\begin{abstract}
We study some classes of lazy cocycles, called pure (respectively neat), 
together with their categorical counterparts, entwined (respectively 
strongly entwined) monoidal categories.   
\end{abstract}
\section*{Introduction}
Let $H$ be a Hopf algebra with bijective antipode. A left 2-cocycle 
$\sigma :H\ot H\rightarrow k$ is called {\it lazy} if it satisfies the 
condition
\begin{eqnarray*}
&&\sigma (h_1, h'_1)h_2h'_2=h_1h'_1\sigma (h_2, h'_2), \;\;\;\;
\forall \; h, h'\in H.
\end{eqnarray*}
In dual form (and with a different name) lazy cocycles appear for instance 
in Majid's book \cite{m}; their most important property, the fact that 
they form a group (denoted now by $Z^2_L(H)$) appears in the 
paper of Chen \cite{chen}. Present terminology stems from   
\cite{bichon}, \cite{cc}, inspired by the fact that a Doi twisting by a 
lazy cocycle does not modify $H$. Moreover, one may define lazy  
2-coboundaries $B^2_L(H)$ and the second lazy cohomology group  
$H^2_L(H)=Z^2_L(H)/B^2_L(H)$, generalizing Sweedler's second 
cohomology group of a cocommutative Hopf algebra (this is done 
by Schauenburg in \cite{sch}).\\
Lazy cocycles have been studied systematically in \cite{bichon},  
\cite{c},  
\cite{cp}, also in connection with Brauer groups of Hopf  
algebras, Bigalois groups, projective representations.\\
In this paper we study a certain class of lazy cocycles,     
satisfying the condition 
\begin{eqnarray*}
&&\sigma (ab_1, c_1)\sigma ^{-1}(b_2, c_2)\sigma (b_3, c_3d)=
\sigma (b_1, c_1d)\sigma ^{-1}(b_2, c_2)\sigma (ab_3, c_3), 
\end{eqnarray*}
for all $a, b, c, d\in H$, called {\it pure} lazy cocycles.  
In dual form, they have been introduced in \cite{doru} as    
{\it pure-braided structure}. This purity condition has a topological 
meaning: pure lazy cocycles give rise to representations of pure  
braid groups and invariants for long knots, cf. \cite{doru}.\\ 
A natural problem is whether it is possible to determine all 
pure lazy cocycles on a given Hopf algebra; 
this seems to be complicated even for ``easy''   
Hopf algebras. This is why we have looked for a stronger condition than  
purity, and we were led to the following concept: a lazy cocycle is called 
{\it neat} if it satisfies the condition 
\begin{eqnarray*}
&&\sigma (a, b_1)\sigma (b_2, c)=\sigma (b_1, c)\sigma (a, b_2), 
\end{eqnarray*}
for all $a, b, c\in H$. It turns out that a neat lazy cocycle is pure, and,   
using the description of lazy cocycles for Sweedler's Hopf algebra 
$H_4$ from \cite{bichon}, it is quite easy to see that any lazy cocycle 
for $H_4$ is neat (hence also pure). \\
The categorical counterpart of pure lazy cocycles was introduced in 
\cite{doru} as a {\it pure-braided category} and  
independently in \cite{brug} as {\it entwined category}. These 
concepts look different but we prove here that they are equivalent (and 
provide another equivalent formulation). We introduce the    
categorical analogue of neat lazy cocycles, as a 
{\it strongly entwined category}, and prove that strongly entwined  
implies entwined. We show that there exists a canonical  
way to produce a strong twine starting from a D-structure  
(consisting of isomorphisms) in the sense of \cite{street}. \\
A natural question is to see what kind of algebraic properties pure and neat  
lazy cocycles have. It turns out that their algebraic properties are not too  
good (for instance they do not seem to form subgroups of   
$Z^2_L(H)$), but are also not so bad, for instance they have a good 
behaviour when extending to a Drinfeld double or a Radford biproduct (this 
extension property may also be regarded as a potential source of examples 
of pure and neat lazy cocycles).     
\section{Preliminaries}\label{sec1}
\setcounter{equation}{0}
In this section we recall basic definitions and results and we fix 
notation to be used throughout the paper. 
All algebras, linear 
spaces, etc, will be over a base field $k$; unadorned $\ot$ means $\ot _k$. 
For a Hopf algebra $H$ with comultiplication $\Delta$ we use 
Sweedler's sigma notation: $\Delta (h)=h_1\ot h_2$ or 
$\Delta (h)=h_{(1)}\ot h_{(2)}$. 
Unless otherwise stated, $H$ will denote a Hopf algebra with
bijective antipode $S$. For a linear map $\sigma :H\ot
H\rightarrow k$ we use either the notation $\sigma (h,h')$
or $\sigma(h\ot h')$. For terminology concerning Hopf algebras and 
monoidal categories we refer to \cite{k}, \cite{m}, \cite{mon}, \cite{sw}. \\
A linear map $\sigma :H\ot H\rightarrow k$ is called a {\it left
2-cocycle} if it satisfies the condition
\begin{eqnarray}
&&\sigma (a_1, b_1)\sigma (a_2b_2, c)=\sigma (b_1, c_1)\sigma (a, 
b_2c_2), \label{leftco}
\end{eqnarray}
for all $a, b, c\in H$, and it is called a {\it right 
2-cocycle} if it satisfies the condition 
\begin{eqnarray}
&&\sigma (a_1b_1, 
c)\sigma (a_2, b_2)=\sigma (a, b_1c_1)\sigma (b_2, c_2). \label{rightco}
\end{eqnarray}
Given 
a linear map $\sigma :H\ot H\rightarrow k$, define a product
$\cdot_{\sigma }$ on $H$ by
\begin{eqnarray*}
&&h\cdot _{\sigma }h'=
\sigma (h_1, h'_1)h_2h'_2, \;\;\;\;\;\forall \;h, h'\in H.
\end{eqnarray*}
Then $\cdot_{\sigma }$ is associative if and only if $\sigma $ is
a left 2-cocycle. If we define $\cdot _{\sigma }$ by
\begin{eqnarray*}
&&h\cdot _{\sigma }h'=h_1h'_1\sigma (h_2, h'_2), \;\;\;\;\;\forall \;
h, h'\in H,
\end{eqnarray*}
then $\cdot _{\sigma }$ is associative if and only if $\sigma $
is a right 2-cocycle. In any of the two cases, $\sigma$ is
normalized (i.e. $\sigma (1, h)= \sigma (h, 1)=\varepsilon (h)$
for all $h\in H$) if and only if $1_H$ is the unit for $\cdot
_{\sigma }$.  If $\sigma $ is a normalized left (respectively
right) 2-cocycle, we denote the algebra $(H, \cdot _{\sigma })$
by $_{\sigma }H$ (respectively $H_{\sigma }$). It is well-known
that $_{\sigma }H$ (respectively $H_{\sigma }$) is a right
(respectively left) $H$-comodule algebra via the comultiplication
$\Delta $ of $H$. If $\sigma :H\ot H\rightarrow k$ is normalized
and convolution invertible, then $\sigma $ is a left 2-cocycle if
and only if $\sigma ^{-1}$ is a right 2-cocycle. \\
If $\gamma :H\rightarrow k$ is linear, normalized (i.e. $\gamma
(1)=1$) and convolution invertible, define
\begin{eqnarray*}
&&D^1(\gamma ):H\ot H\rightarrow k, \;\;\;
D^1(\gamma )(h, h')=\gamma (h_1)\gamma (h'_1)\gamma ^{-1}(h_2h'_2),
\;\;\;\;\;\forall \;h, h'\in H.
\end{eqnarray*}
Then $D^1(\gamma )$ is a normalized and convolution invertible
left 2-cocycle.\\
We recall from \cite{bichon} some facts about lazy cocycles and  
lazy cohomology. The set $Reg^1 (H)$ (respectively $Reg^2 (H)$)
consisting of normalized and convolution invertible linear maps
$\gamma :H\rightarrow k$ (respectively $\sigma :H\ot H\rightarrow
k$), is a group with respect to the convolution product. 
An element $\gamma 
\in Reg^1 (H)$ is called {\it lazy} if
\begin{eqnarray}
&&\gamma (h_1)h_2=h_1\gamma (h_2), \;\;\;\;\;\forall \;h\in H. \label{lazy1}
\end{eqnarray}
The set of lazy elements of $Reg^1 (H)$, denoted by $Reg^1_L (H)$,
is a central subgroup of $Reg^1 (H)$. An element $\sigma \in
Reg^2 (H)$ is called {\it lazy} if
\begin{eqnarray}
&&\sigma (h_1, h'_1)h_2h'_2=h_1h'_1\sigma (h_2, h'_2),\;\;\;\;\;
\forall \;h, h'\in H. \label{lazy2}
\end{eqnarray}
The set of lazy elements of $Reg^2 (H)$, denoted by $Reg^2_L (H)$,
is a subgroup of $Reg^2 (H)$. We denote by $Z^2 (H)$ the set of
left 2-cocycles on $H$ and by $Z^2_L (H)$ the set $Z^2 (H)\cap
Reg^2_L (H)$ of normalized and convolution invertible lazy
2-cocycles. If $\sigma \in Z^2_L(H)$, then the algebras $_{\sigma
}H$ and $H_{\sigma }$ coincide and will be denoted by $H(\sigma
)$; moreover, $H(\sigma )$ is an $H$-bicomodule algebra via
$\Delta $. \\
It is well-known that in general the set $Z^2 (H)$ of left
2-cocycles is not closed under convolution. One of the main
features of lazy 2-cocycles is that the set $Z^2_L (H)$ is closed
under convolution, and that the convolution inverse of an element
$\sigma \in Z^2_L (H)$ is again a lazy 2-cocycle, so $Z^2_L (H)$
is a group under convolution. In particular, a lazy 2-cocycle is
also a right 2-cocycle.\\ 
Consider now the map 
$D^1:Reg^1 (H)\rightarrow Reg^2 (H)$, $D^1(\gamma )(h, h')=\gamma  
(h_1)\gamma (h'_1)\gamma ^{-1}(h_2h'_2)$, 
for all $h, h'\in H$. 
Then, by \cite{bichon}, the map $D^1$ induces a group morphism 
$Reg^1_L (H)\rightarrow Z^2_L (H)$, with image contained in 
the centre of $Z^2_L (H)$; denote by $B^2_L (H)$ this central
subgroup $D^1(Reg^1_L (H))$ of $Z^2_L (H)$ (its elements are
called {\it lazy 2-coboundaries}). Then define the {\it second 
lazy cohomology group} $H^2_L (H)=Z^2_L (H)/B^2_L (H)$.
\section{Pure-braided and entwined monoidal categories}
\setcounter{equation}{0}
We begin this section by recalling the following two concepts (all 
monoidal categories are assumed to be strict, with unit denoted by $I$).
\begin{definition} (\cite{doru}) Let ${\cal C}$  be a monoidal category.
A pure-braided structure of ${\cal C}$ consists of two families of natural  
isomorphisms $A_{U,V,W}:U\otimes V\otimes W \to U\otimes V\otimes W$ 
and $B_{U,V,W}:U\otimes V\otimes W \to 
U\otimes V\otimes W$ such that:
\begin{eqnarray}
&&A_{U\otimes V,W,X}=A_{U,V\otimes W,X}(id_U\otimes A_{V,W,X}),  
\label{a1} \\
&&A_{U,V,W\otimes X}=(A_{U,V,W}\otimes id_X) A_{U,V\otimes W,X}, 
\label{a2} \\
&&B_{U\otimes V,W,X}=(id_U \otimes B_{V,W,X})B_{U,V\otimes W,X}, 
\label{baba} \\
&&B_{U,V,W\otimes X}=B_{U,V\otimes W,X}(B_{U,V,W}\otimes id_X),  
\label{b2} \\
&&(A_{U,V,W}\otimes id_X)(id_U \otimes B_{V,W,X})=
(id_U \otimes B_{V,W,X})(A_{U,V,W}\otimes id_X), 
\label{cab} \\
&&A_{U,I,V}=B_{U,I,V}.  
\label{t1t} 
\end{eqnarray}
A category equipped with a pure-braided structure is called a pure-braided  
category. 
\end{definition}
\begin{remark} The axioms (\ref{a1})-(\ref{b2}) imply also the 
following relations:
\begin{eqnarray}
&&A_{I,U,V}=A_{U,V,I}=id_{U\otimes V},  
\label{t1a} \\
&&B_{I,U,V}=B_{U,V,I}=id_{U\otimes V}.  
\label{t1b}
\end{eqnarray}
\end{remark}
\begin{definition} (\cite{brug}) Let ${\cal C}$ be  a monoidal category.  
A twine of ${\cal C}$ is a natural isomorphism 
$D_{X,Y}:X\otimes Y\to X\otimes Y$ 
satisfying the following axioms:
\begin{eqnarray*}
&&D_{I,I}=id_I, \label{db0}\\
&&(D_{X,Y}\otimes id_Z)D_{X\otimes Y, Z}=
(id_X \otimes D_{Y,Z})D_{X,Y\otimes Z}, 
\label{db1}\\
&&(D_{X\otimes Y,Z}\otimes id_T)
(id_X\otimes D_{Y,Z}^{-1}\otimes id_T)(id_X\otimes D_{Y,Z\otimes T})\\
&&\;\;\;\;\;\;\;\;\;\;=(id_X\otimes D_{Y,Z\otimes T})
(id_X\otimes D_{Y,Z}^{-1}\otimes id_T)(D_{X\otimes Y,Z}\otimes id_T). 
\label{db2}
\end{eqnarray*}
\label{ent}
A category equipped with a twine is called an entwined category. 
\end{definition}
\begin{remark} By \cite{brug}, if  
$({\cal C}, D)$ is an entwined category then  
$D_{X,I}=D_{I,X}=id_X$, $\;\forall \; X\in {\cal C}$. \label{re1}
\end{remark}
\begin{remark}
If ${\cal C}$ is a monoidal category and $D_{X, Y}:X\ot Y\rightarrow X\ot Y$ 
is a natural isomorphism, the naturality of $D$ implies (for all 
$X, Y, Z\in {\cal C}$): 
\begin{eqnarray}
&&(D_{X, Y}\ot id_Z)D_{X\ot Y, Z}=D_{X\ot Y, Z}(D_{X, Y}\ot id_Z), 
\label{lac1} \\
&&(id_X\ot D_{Y, Z})D_{X, Y\ot Z}=D_{X, Y\ot Z}(id_X\ot D_{Y, Z}). 
\label{lac2}
\end{eqnarray}
\end{remark}
We prove now that these two concepts are equivalent.
\begin{proposition} Let ${\cal C}$ be a monoidal category.  \\
a) If $({\cal C}, A, B)$ is a pure-braided category and we define 
$D_{U, V}:U\otimes V\rightarrow U\otimes V$ by $D_{U, V}:=A_{U, I, V}=
B_{U, I, V}$, then $D_{U, V}$ is a natural isomorphism satisfying 
\begin{eqnarray}
&&D_{I,X}=D_{X,I}=id_X, \label{inter1}
\end{eqnarray}
\begin{eqnarray*}
&&(D_{X, Y}\otimes id_Z\otimes id_T)
(id_X\otimes D_{Y\otimes Z,T})(D_{X\otimes Y, Z}\otimes id_T)
\end{eqnarray*}
\begin{eqnarray}
&&\;\;\;\;\;\;\;\;=(id_X\otimes id_Y \otimes D_{Z, T})(D_{X,Y\otimes Z}
\otimes id_T)(id_X \otimes D_{Y,Z\otimes T}).\label{inter2}
\end{eqnarray} 
b) If $D_{U, V}:U\otimes V\rightarrow U\otimes V$ is a natural isomorphism 
satisfying (\ref{inter1}) and (\ref{inter2}), then $({\cal C},D)$ is 
an entwined category.\\
c) If $({\cal C}, D)$ is an entwined category and we define 
$A_{X, Y, Z}, B_{X, Y, Z}:X\otimes Y\otimes Z\rightarrow 
X\otimes Y\otimes Z$ by 
\begin{eqnarray} 
&&A_{X,Y,Z}=D_{X\otimes Y, Z}(id_X \otimes D_{Y,Z}^{-1})=
(D_{X, Y}^{-1}\otimes id_Z)D_{X,Y\otimes Z}, \label{auri}\\
&&B_{X,Y,Z}=(id_X \otimes D_{Y,Z}^{-1})D_{X\otimes Y, Z}=
D_{X,Y\otimes Z}(D_{X, Y}^{-1}\otimes id_Z), \label{buri}
\end{eqnarray}
then $({\cal C}, A, B)$ is a pure-braided category. 
\label{abd}
\end{proposition} 
\begin{proof}
a) Define $D_{U,V}:=A_{U,I,V}=B_{U, I, V}$. By   
(\ref{a1}) we have 
$A_{U\otimes I,I,X}=A_{U,I\otimes I,X}(id_U\otimes A_{I,I,X})$, hence  
we obtain $D_{I,X}=A_{I,I,X}=id_X$   
and similarly $D_{X,I}=id_X$. We prove that 
$A_{U, V, X}=D_{U\otimes V, X}(id_U\otimes D_{V,X}^{-1})$; 
indeed, we have:
\begin{eqnarray*}
D_{U\otimes V, X}&=&A_{U\otimes V,I,X}\\
{\rm (\ref{a1})}&=&A_{U,V\otimes I,X}(id_U\otimes A_{V,I,X})\\
&=&A_{U, V, X}(id_U\otimes D_{V,X}), 
\end{eqnarray*}
and similarly
\begin{eqnarray*}
&&A_{U, V, X}=(D_{U,V}^{-1}\otimes id_X)D_{U, V\otimes  X}, \\
&&B_{U, V, X}=(id_U\otimes D_{V,X}^{-1})D_{U\otimes V, X}, \\
&&B_{U, V, X}=D_{U, V\otimes  X}(D_{U,V}^{-1}\otimes id_X).
\end{eqnarray*} 
Using these formulae we obtain:
\begin{eqnarray*}
&&(A_{U,V,W}\otimes id_X)(id_U \otimes B_{V,W,X})\\
&&\;\;\;\;\;
=(D_{U,V}^{-1}\otimes id_W\otimes id_X)(D_{U, V\otimes  W}\otimes id_X)
(id_U\otimes D_{V, W\otimes X})(id_U\otimes D_{V,W}^{-1}\otimes id_X), \\
&&(id_U \otimes B_{V,W,X})(A_{U,V,W}\otimes id_X)\\
&&\;\;\;\;\;=(id_U\otimes id_V \otimes  D_{W,X}^{-1})
(id_U\otimes D_{V\otimes W, X})(D_{U\otimes V, W}\otimes id_X) 
(id_U\otimes D_{V,W}^{-1}\otimes id_X).
\end{eqnarray*}
Now using (\ref{cab}) we get (\ref{inter2}).\\
b) We take $T=I$  in (\ref{inter2}), obtaining    
\begin{eqnarray*}
&&(D_{X, Y}\otimes id_Z\otimes id_I)
(id_X\otimes D_{Y\otimes Z,I})(D_{X\otimes Y, Z}\otimes id_I)\\
&&\;\;\;\;\;=(id_X\otimes id_Y \otimes D_{Z, I})(D_{X,Y\otimes Z}\otimes id_I)
(id_X \otimes D_{Y,Z\otimes I}), 
\end{eqnarray*}
which can be rewritten as
\begin{eqnarray}
(D_{X, Y}\otimes id_Z)D_{X\otimes Y, Z}=
D_{X,Y\otimes Z}(id_X \otimes D_{Y,Z}).
\label{123}
\end{eqnarray}
Also, (\ref{inter2}) implies    
\begin{eqnarray*}
&&(id_X\otimes id_Y \otimes D_{Z, T}^{-1})
(id_X\otimes D_{Y\otimes Z,T})(D_{X\otimes Y, Z}\otimes id_T)\\
&&\;\;\;\;\;=(D_{X, Y}^{-1}\otimes id_Z\otimes id_T)
(D_{X,Y\otimes Z}\otimes id_T)(id_X \otimes D_{Y,Z\otimes T}), 
\end{eqnarray*}
and using (\ref{123}) we obtain
\begin{eqnarray*}
&&(id_X\otimes D_{Y,Z\otimes T})(id_X\otimes D_{Y, Z}^{-1}\otimes id_T)
(D_{X\otimes Y, Z}\otimes id_T)\\
&&\;\;\;\;\;=(D_{X\otimes Y, Z}\otimes id_T)
(id_X\otimes D_{Y, Z}^{-1}\otimes id_T)
(id_X \otimes D_{Y,Z\otimes T}).
\end{eqnarray*}
c) Define $A$ and $B$ by (\ref{auri}) and (\ref{buri}) respectively. 
We prove (\ref{a1}):
\begin{eqnarray*}
A_{U\otimes V,W,X}&=&D_{U\otimes V\otimes W, X}
(id_{U\otimes V}\otimes D_{W,X}^{-1})\\
&=&D_{U\otimes V\otimes W, X}(id_U\otimes D_{V\otimes W,X}^{-1})
(id_U\otimes D_{V\otimes W,X})(id_{U\otimes V}\otimes D_{W,X}^{-1})\\
&=&A_{U,V\otimes W,X}(id_U\otimes A_{V,W,X}).
\end{eqnarray*}
Similarly we get (\ref{a2}), (\ref{baba}) and (\ref{b2}).  
From the definition we have $A_{U,I,V}=D_{U,V}=B_{U,I,V}$. 
Finally, we prove (\ref{cab}): \\[2mm]
${\;\;\;\;\;}$
$(A_{U,V,W}\otimes id_X)(id_U \otimes B_{V,W,X})$
\begin{eqnarray*}
&&=(D_{U\otimes V, W}\otimes id_X)(id_U \otimes D_{V,W}^{-1}\otimes id_X) 
(id_U\otimes D_{V, W\otimes X})(id_U \otimes D_{V,W}^{-1}\otimes id_X)\\
&&=(id_U\otimes D_{V, W\otimes X})(id_U \otimes D_{V,W}^{-1}\otimes id_X)
(D_{U\otimes V, W}\otimes id_X)(id_U \otimes D_{V,W}^{-1}\otimes id_X)\\
&&=(id_U \otimes B_{V,W,X})(A_{U,V,W}\otimes id_X),
\end{eqnarray*}
finishing the proof.
\end{proof} 
\begin{definition} Let ${\cal C}$ be a monoidal category 
and $T_{U, V}:U\otimes V\to U\otimes V$ a natural isomorphism.  
We say that $T$ is a strong twine (or $({\cal C}, T)$ is 
strongly entwined) if:
\begin{eqnarray}
&&T_{I,I}=id_I, \label{st1} \\
&&(T_{U, V}\otimes id_W)T_{U\otimes V, W}=
(id_U\otimes T_{V, W})T_{U, V\otimes W}, \label{st2} \\
&&(T_{U, V}\otimes id_W)(id_U\otimes T_{V, W})=
(id_U\otimes T_{V, W})(T_{U, V}\otimes id_W). \label{st3}
\end{eqnarray}
\end{definition}
\begin{proposition}\label{lilu} If $({\cal C}, T)$ is strongly entwined then  
$({\cal C}, T)$ is entwined.
\end{proposition}
\begin{proof}  
First we prove that 
\begin{eqnarray} 
&&(T_{U,V}\otimes id_{W\otimes X})(id_U\otimes T_{V\otimes W, X})=
(id_U\otimes T_{V\otimes W, X})(T_{U,V}\otimes id_{W\otimes X}).
\label{con1}
\end{eqnarray}
Indeed, we have:\\[2mm]
${\;\;\;\;\;}$
$(T_{U,V}\otimes id_{W\otimes X})(id_U\otimes T_{V\otimes W, X})$
\begin{eqnarray*}
{\rm (\ref{st2})}&=&(T_{U,V}\otimes id_{W\otimes X})
(id_U\otimes T_{V,W}^{-1}\otimes id_X)
(id_{U\otimes V}\otimes T_{W,X})(id_U\otimes T_{V,W\otimes X})\\
{\rm (\ref{st3})}&=&(id_U\otimes T_{V,W}^{-1}\otimes id_X)
(T_{U,V}\otimes id_{W\otimes X})
(id_{U\otimes V}\ot T_{W, X})(id_U\otimes T_{V, W\ot X})\\
&=&(id_U\otimes T_{V,W}^{-1}\otimes id_X)(id_{U\otimes V}\ot 
T_{W, X})
(T_{U,V}\otimes id_{W\otimes X})(id_U\otimes T_{V, W\ot X})\\
{\rm (\ref{st3})}&=&(id_U\otimes T_{V,W}^{-1}\otimes id_X)
(id_{U\otimes V}\ot T_{W, X})
(id_U\otimes T_{V, W\ot X})(T_{U,V}\otimes id_{W\otimes X})\\
{\rm (\ref{st2})}&=&(id_U\otimes T_{V\otimes W, X})
(T_{U,V}\otimes id_{W\otimes X}), 
\end{eqnarray*}
and similarly
\begin{eqnarray}
&&(T_{U,V\otimes W}\otimes id_X)(id_{U\otimes V}\otimes T_{W, X})=
(id_{U\otimes V}\otimes T_{W, X})(T_{U,V\otimes W}\otimes id_X).
\label{con2}
\end{eqnarray}
Now we compute:\\[2mm]
${\;\;\;\;}$
$(T_{U\otimes V,W}\otimes id_X)(id_U\otimes T_{V, W}^{-1}\otimes id_X)
(id_U\otimes T_{V,W\otimes X})$
\begin{eqnarray*}
{\rm (\ref{st2})}&=&(T_{U,V}^{-1}\otimes id_{W\otimes X})
(id_U\otimes T_{V, W}\otimes id_X)(T_{U,V\otimes W}\otimes id_X)
(id_U\otimes T_{V, W}^{-1}\otimes id_X)\\
&&(id_{U\otimes V}\otimes T_{W, X}^{-1})(id_U\otimes T_{V, W}\otimes id_X)
(id_U\otimes T_{V\otimes W,X})\\
{\rm (\ref{lac2})}&=&(T_{U,V}^{-1}\otimes id_{W\otimes X})
(T_{U,V\otimes W}\otimes id_X)
(id_{U\otimes V}\otimes T_{W, X}^{-1})\\
&&(id_U\otimes T_{V\otimes W,X})
(id_U\otimes T_{V, W}\otimes id_X)\\
{\rm (\ref{con2})}&=&(T_{U,V}^{-1}\otimes id_{W\otimes X})(id_{U\otimes V}
\otimes T_{W, X}^{-1})(T_{U,V\otimes W}\otimes id_X)\\
&&(id_U\otimes T_{V\otimes W,X})(id_U\otimes T_{V, W}\otimes id_X)\\
{\rm (\ref{st3})}&=&(T_{U,V}^{-1}\otimes id_{W\otimes X})(id_{U\otimes V}
\otimes T_{W, X}^{-1})(id_U\otimes T_{V\otimes W,X})\\
&&(T_{U,V\otimes W}\otimes id_X)(id_U\otimes T_{V, W}\otimes id_X), 
\end{eqnarray*}
${\;\;\;\;}$
$(id_U\otimes T_{V,W\otimes X})(id_U\otimes T_{V, W}^{-1}\otimes id_X)
(T_{U\otimes V,W}\otimes id_X)$
\begin{eqnarray*}
{\rm (\ref{st2})}&=&(id_{U\otimes V}\otimes T_{W, X}^{-1})
(id_U\otimes T_{V, W}
\otimes id_X)(id_U\otimes T_{V\otimes W,X})
(id_U\otimes T_{V, W}^{-1}\otimes id_X)\\
&&(T_{U,V}^{-1}\otimes id_{W\otimes X})(id_U\otimes T_{V, W}\otimes id_X)
(T_{U,V\otimes W}\otimes id_X)\\
{\rm (\ref{lac1})}&=&(id_{U\otimes V}\otimes T_{W, X}^{-1})
(id_U\otimes T_{V\otimes W,X})
(T_{U,V}^{-1}\otimes id_{W\otimes X})\\
&&(id_U\otimes T_{V, W}\otimes id_X)(T_{U,V\otimes W}\otimes id_X)\\
{\rm (\ref{con1}, \ref{lac2})}&=&(id_{U\otimes V}\otimes T_{W, X}^{-1})
(T_{U,V}^{-1}\otimes id_{W\otimes X})(id_U\otimes T_{V\otimes W,X})\\
&&(T_{U,V\otimes W}\otimes id_X)(id_U\otimes T_{V, W}\otimes id_X), 
\end{eqnarray*}
showing that $({\cal C}, T)$ is an entwined category.
\end{proof}
\begin{remark}
Any monoidal category contains at least one strong twine: the 
trivial one.
\end{remark}
The categorical analogue of the operator $D^1$ from the Preliminaries 
looks as follows (see \cite{brug}). 
If ${\cal C}$ is a monoidal category and   
$R_X:X\rightarrow X$ is a natural isomorphism in ${\cal C}$ such that  
$R_I=id_I$, we denote $D^1(R)_{X, Y}:=(R_X\ot R_Y)R_{X\ot Y}^{-1}= 
R_{X\ot Y}^{-1}(R_X\ot R_Y)$ as a morphism $X\otimes Y
\rightarrow X\otimes Y$. 
\begin{definition} (\cite{street})  
Let ${\cal C}$ be a monoidal category. A D-structure on 
${\cal C}$ consists of a family of natural morphisms $R_X:X\rightarrow X$  
in ${\cal C}$, such that $R_I=id_I$ and (for all $X, Y, Z\in {\cal C}$):
\begin{eqnarray}
&&(R_{X\ot Y}\ot id_Z)(id_X\ot R_{Y\ot Z})=(id_X\ot R_{Y\ot Z})
(R_{X\ot Y}\ot id_Z). \label{dstr}
\end{eqnarray}
\end{definition}
This concept provides a method for constructing strong twines, as follows:
\begin{proposition} Let ${\cal C}$ be a monoidal category and $R$ a 
D-structure on ${\cal C}$ such that all $R_X$ are isomorphisms. Then 
$D^1(R)$ is a strong twine on ${\cal C}$.
\end{proposition} 
\begin{proof}
We only have to check (\ref{st3}). We compute: \\[2mm]
${\;\;\;\;}$
$(D^1(R)_{U, V}\ot id_W)(id_U\ot D^1(R)_{V, W})$
\begin{eqnarray*}
&=&((R_U\ot R_V)R_{U\ot V}^{-1}\ot id_W)(id_U\ot R_{V\ot W}^{-1}(R_V\ot 
R_W))\\
&=&(R_U\ot R_V\ot id_W)(R_{U\ot V}^{-1}\ot id_W)(id_U\ot R_{V\ot W}^{-1})
(id_U\ot R_V\ot R_W)\\
{\rm (\ref{dstr})}&=&(R_U\ot R_V\ot id_W)(id_U\ot R_{V\ot W}^{-1})
(R_{U\ot V}^{-1}\ot id_W)(id_U\ot R_V\ot R_W)\\
&=&(R_U\ot (R_V\ot id_W)R_{V\ot W}^{-1})(R_{U\ot V}^{-1}
(id_U\ot R_V)\ot R_W)\\
&=&(R_U\ot R_{V\ot W}^{-1}(R_V\ot id_W))((id_U\ot R_V)R_{U\ot V}^{-1}\ot 
R_W)\\
&=&(id_U\ot R_{V\ot W}^{-1})(R_U\ot R_V\ot id_W)(id_U\ot R_V\ot R_W)
(R_{U\ot V}^{-1}\ot id_W)\\
&=&(id_U\ot R_{V\ot W}^{-1})(id_U\ot R_V\ot R_W)(R_U\ot R_V\ot id_W)
(R_{U\ot V}^{-1}\ot id_W)\\
&=&(id_U\ot D^1(R)_{V, W})(D^1(R)_{U, V}\ot id_W), 
\end{eqnarray*}
finishing the proof. 
\end{proof} 
\section{Pure and neat lazy cocycles}
\setcounter{equation}{0}
\begin{definition}{\rm
Let $\sigma \in Reg^2(H)$; we call $\sigma $ {\it pure} if it 
satisfies the condition:
\begin{eqnarray}
&&\sigma (ab_1, c_1)\sigma ^{-1}(b_2, c_2)\sigma (b_3, c_3d)=
\sigma (b_1, c_1d)\sigma ^{-1}(b_2, c_2)\sigma (ab_3, c_3), \label{pure1}
\end{eqnarray}
for all $a, b, c\in H$. If $\sigma $ is moreover lazy we call it pure lazy 
and denote by $Reg ^2_{PL}(H)$ the set of pure lazy elements. We also 
denote by $Z^2_{PL}(H)$ the set of pure lazy 2-cocycles.}
\end{definition}
\begin{remark}{\em
The concept of pure lazy cocycle is dual to the concept of 
{\it pure-braided structure} in \cite{doru}.}
\end{remark}
\begin{example}{\em 
If $r, s$ are two coquasitriangular structures on $H$, then $r_{21}*s$   
is a pure lazy 2-cocycle. The fact that it is a lazy 2-cocycle was noticed 
in \cite{bichon}, and the fact that it is pure is analogous to a remark 
due to Virelizier, see \cite{brug}.}
\end{example} 
\begin{definition}{\rm
Let $\sigma \in Reg^2(H)$; we call $\sigma $ {\it neat} if it   
satisfies the condition:
\begin{eqnarray}
&&\sigma (a, b_1)\sigma (b_2, c)=
\sigma (b_1, c)\sigma (a, b_2), \label{verypure1}
\end{eqnarray}
for all $a, b, c\in H$. If $\sigma $ is moreover lazy we call it 
neat lazy   
and denote by $Reg ^2_{NL}(H)$ the set of neat lazy elements. We also   
denote by $Z^2_{NL}(H)$ the set of neat lazy 2-cocycles.} 
\end{definition}
\begin{remark}{\em
Relation (\ref{verypure1}) is a commutation condition. Namely, define the 
maps $\varphi , \psi :H\rightarrow H^*$, $\varphi (a)(b)=\sigma (a, b)$ 
and $\psi (a)(b)=\sigma (b, a)$. Then (\ref{verypure1}) holds if and 
only if $\varphi (a)*\psi (c)=\psi (c)*\varphi (a)$ in $H^*$, for all 
$a, c\in H$.} 
\end{remark}
We have the following dictionary between lazy cocycles and categorical 
structures:
\begin{proposition} Let $H$ be a Hopf algebra, $\sigma \in Reg^2_L(H)$, 
and consider ${\cal C}={\cal M}^H$, the category of right $H$-comodules, 
with tensor product given by $(m\ot n)_{(0)}\ot (m\ot n)_{(1)}=
(m_{(0)}\ot n_{(0)})\ot m_{(1)}n_{(1)}$. Define 
$T_{M, N}(m\ot n)=m_{(0)}\ot n_{(0)}\sigma (m_{(1)}, n_{(1)})$. Then  
$\sigma $ is a pure (respectively neat) lazy  
2-cocycle if and only if $T$ is a twine (respectively strong twine). 
\end{proposition}
As a consequence of this and Proposition \ref{lilu}, we obtain: 
\begin{proposition}\label{very}
$Z^2_{NL}(H)\subseteq Z^2_{PL}(H)$, that is any neat lazy cocycle is pure. 
\end{proposition}
\begin{remark}{\em 
A pure lazy cocycle of the type $r_{21}*s$,  
with $r, s$ coquasitriangular  
structures on $H$, is not necessarily neat.}
\end{remark} 
\begin{example}{\em 
Let $H_4$ be Sweedler's 4-dimensional Hopf algebra. A description of 
$Z^2_L(H_4)$ was given in \cite{bichon}, Example 2.1. Using the 
formulae in \cite{bichon}, one can prove, by a direct computation, that 
any lazy 2-cocycle on $H_4$ is neat, hence we obtain  
$Z^2_{NL}(H_4)=Z^2_{PL}(H_4)=Z^2_L(H_4)$.}
\end{example}
\begin{proposition} (see \cite{brug}, \cite{doru}) 
If $\gamma \in Reg^1_L(H)$ satisfies the condition \\[2mm]
${\;\;\;\;\;\;\;\;\;\;}$
$\gamma (a_1b_1)\gamma ^{-1}(a_2b_2c_1)\gamma (b_3c_2)
\gamma ^{-1}(b_4c_3d_1)\gamma (c_4d_2)=$
\begin{eqnarray}
&&=\gamma (c_1d_1)\gamma ^{-1}(b_1c_2d_2) 
\gamma (b_2c_3)\gamma ^{-1}(a_1b_3c_4)\gamma (a_2b_4), \label{pure0}
\end{eqnarray}
for all $a, b, c, d\in H$, then $D^1(\gamma )\in Z^2_{PL}(H)$. An element 
$\gamma \in Reg^1(H)$ satisfying (\ref{pure0}) is said to be pure.    
We denote 
by $Reg^1_{PL}(H)$ the set of pure lazy elements. 
\end{proposition}
\begin{definition}
An element $\gamma \in Reg^1(H)$ satisfying the condition  
\begin{eqnarray}
&&\gamma (ab_1)\gamma (b_2c)=\gamma (b_1c)\gamma (ab_2), \label{pure2}
\end{eqnarray}
for all $a, b, c\in H$, is said to be neat. We denote the set 
of neat elements by $Reg^1_{N}(H)$    
and the set of neat lazy elements by $Reg^1_{NL}(H)$ .   
\end{definition}
\begin{remark}
A neat lazy element corresponds to a $D$-morphism in   
\cite{street}, except for the fact that a $D$-morphism is not required to be  
convolution invertible.
\end{remark}
\begin{proposition} \label{vpln}
$Reg^1_{NL}(H)\subseteq Reg^1_{PL}(H)$ and $D^1(Reg^1_{NL}(H))\subseteq  
Z^2_{NL}(H)$.
\end{proposition}
\begin{proof} Straightforward computation.
\end{proof}
\begin{proposition}
If $\gamma \in Reg^1(H)$ satisfies the condition 
\begin{eqnarray}
&&\gamma (ab_1)b_2=\gamma (ab_2)b_1, \label{pure3}
\end{eqnarray}
for all $a, b \in H$, then $\gamma \in Reg^1_{NL}(H)$.
\end{proposition}
\begin{proof}
An element $\gamma $  
satisfying (\ref{pure3}) is automatically lazy and also satisfies   
(\ref{pure2}). 
\end{proof}
\section{Extending pure and neat lazy cocycles to Drinfeld doubles and 
Radford biproducts}
\setcounter{equation}{0}
Let $H$ be a finite dimensional Hopf algebra. 
Recall that the Drinfeld double $D(H)$ is a quasitriangular Hopf 
algebra realized on the $k$-linear space $H^*\ot H$; its
coalgebra structure is $H^{*cop}\ot H$ and the algebra
structure is given by 
\begin{eqnarray*}
(p\ot h)(q\ot l)=p(h_1\rightharpoonup 
q\leftharpoonup S^{-1}(h_3))\ot h_2l,
\end{eqnarray*}
for all $p, q\in H^*$ and 
$h, l\in H$, where $\rightharpoonup $ and $\leftharpoonup $ are
the left and right regular actions of $H$ on $H^*$ given by
$(h\rightharpoonup p)(l)=p(lh)$ and $(p\leftharpoonup
h)(l)=p(hl)$ for all $h, l\in H$ and $p\in H^*$.\\
If $\sigma \in Z^2_L(H)$, define 
$\overline{\sigma }:D(H)\ot D(H)\rightarrow k$ by  
\begin{eqnarray}
&&\overline{\sigma }(p\ot h, q\ot l)=p(1)q(S^{-1}(h_3)h_1)\sigma (h_2, l),
 \label{exti} 
\end{eqnarray}
for all $p, q\in H^*$ and $h, l\in H$. Then, by \cite{cp}, 
$\overline{\sigma }\in Z^2_L(D(H))$, and its convolution inverse is  
\begin{eqnarray}
&&\overline{\sigma }^{-1}(p\ot h, q\ot l)=
p(1)q(S^{-1}(h_3)h_1)\sigma ^{-1}(h_2, l). \label{extiinv} 
\end{eqnarray}
Moreover, we have: 
\begin{proposition}
If $\sigma \in Z^2_{PL}(H)$ then $\overline{\sigma }\in Z^2_{PL}(D(H))$.  
If $\sigma \in Z^2_{NL}(H)$ then $\overline{\sigma }\in Z^2_{NL}(D(H))$.  
\end{proposition}
\begin{proof}
Assume first that $\sigma \in Z^2_{PL}(H)$ and let $a, b, c, d\in H$ and 
$A, B, C, D\in H^*$; we prove (\ref{pure1}) for $\overline{\sigma }$ and 
the elements $A\ot a$, $B\ot b$, $C\ot c$, $D\ot d$ in $D(H)$. 
We compute:\\[2mm]
${\;\;\;\;\;}$$\overline{\sigma }((A\ot a)(B\ot b)_1, (C\ot c)_1)
\overline{\sigma }^{-1}((B\ot b)_2, (C\ot c)_2)
\overline{\sigma }((B\ot b)_3, (C\ot c)_3(D\ot d))$
\begin{eqnarray*}
&=&\overline{\sigma }((A\ot a)(B_3\ot b_1), C_3\ot c_1)
\overline{\sigma }^{-1}(B_2\ot b_2, C_2\ot c_2)\\
&&\overline{\sigma }(B_1\ot b_3, (C_1\ot c_3)(D\ot d))\\
&=&\overline{\sigma }(A(a_1\rightharpoonup B_3
\leftharpoonup S^{-1}(a_3))\ot a_2b_1, C_3\ot c_1)
\overline{\sigma }^{-1}(B_2\ot b_2, C_2\ot c_2)\\
&&\overline{\sigma }(B_1\ot b_3, C_1(c_3\rightharpoonup D
\leftharpoonup S^{-1}(c_5))\ot c_4d)\\
{\rm (\ref{exti}, \ref{extiinv})}&=&A(1)B_3(S^{-1}(a_3)a_1)
C_3(S^{-1}(a_{(2, 3)}b_{(1, 3)})a_{(2, 1)}b_{(1, 1)})
\sigma (a_{(2, 2)}b_{(1, 2)}, c_1)\\
&&B_2(1)C_2(S^{-1}(b_{(2, 3)})
b_{(2, 1)})\sigma ^{-1}(b_{(2, 2)}, c_2)\\
&&B_1(1)C_1(S^{-1}(b_{(3, 3)})_1b_{(3, 1, 1)})
D(S^{-1}(c_5)S^{-1}(b_{(3, 3)})_2b_{(3, 1, 2)}c_3)
\sigma (b_{(3, 2)}, c_4d)\\
&=&A(1)B(S^{-1}(a_5)a_1)C(S^{-1}(b_{11})b_7S^{-1}(b_6)b_4
S^{-1}(b_3)S^{-1}(a_4)a_2b_1)\\
&&D(S^{-1}(c_5)S^{-1}(b_{10})b_8c_3)
\sigma (a_3b_2, c_1)\sigma ^{-1}(b_5, c_2)\sigma (b_9, c_4d)\\
&=&A(1)B(S^{-1}(a_5)a_1)C(S^{-1}(b_7)S^{-1}(a_4)a_2b_1)
D(S^{-1}(b_6c_5)b_4c_3)\\
&&\sigma (a_3b_2, c_1)\sigma ^{-1}(b_3, c_2)\sigma (b_5, c_4d)\\
{\rm (\ref{lazy2})}&=&A(1)B(S^{-1}(a_5)a_1)C(S^{-1}(b_7)S^{-1}(a_4)a_2b_1)
D(S^{-1}(b_6c_5)b_5c_4d_3S^{-1}(d_1))\\
&&\sigma (a_3b_2, c_1)\sigma ^{-1}(b_3, c_2)\sigma (b_4, c_3d_2)\\
{\rm (\ref{pure1})}&=&A(1)B(S^{-1}(a_5)a_1)C(S^{-1}(b_5)S^{-1}(a_4)a_2b_1)
D(d_3S^{-1}(d_1))\\
&&\sigma (b_2, c_1d_2)\sigma ^{-1}(b_3, c_2)\sigma (a_3b_4, c_3), 
\end{eqnarray*}
${\;\;\;\;\;}$$\overline{\sigma }((B\ot b)_1, (C\ot c)_1(D\ot d))
\overline{\sigma }^{-1}((B\ot b)_2, (C\ot c)_2)\overline{\sigma }((A\ot a)
(B\ot b)_3, (C\ot c)_3)$
\begin{eqnarray*}
&=&\overline{\sigma }(B_3\ot b_1, (C_3\ot c_1)(D\ot d))
\overline{\sigma }^{-1}(B_2\ot b_2, C_2\ot c_2)\\
&&\overline{\sigma }((A\ot a)(B_1\ot b_3), C_1\ot c_3)\\
&=&\overline{\sigma }(B_3\ot b_1, C_3(c_1\rightharpoonup D
\leftharpoonup S^{-1}(c_3))\ot c_2d)
\overline{\sigma }^{-1}(B_2\ot b_2, C_2\ot c_4)\\
&&\overline{\sigma }(A(a_1\rightharpoonup B_1
\leftharpoonup S^{-1}(a_3))\ot a_2b_3, C_1\ot c_5)\\
{\rm (\ref{exti}, \ref{extiinv})}&=&B_3(1)C_3(S^{-1}(b_{(1, 3)})_1
b_{(1, 1, 1)})D(S^{-1}(c_3)S^{-1}(b_{(1, 3)})_2b_{(1, 1, 2)}c_1)
\sigma (b_{(1, 2)}, c_2d)\\
&&B_2(1)C_2(S^{-1}(b_{(2, 3)})b_{(2, 1)})\sigma ^{-1}(b_{(2, 2)}, c_4)\\
&&A(1)B_1(S^{-1}(a_3)a_1)C_1(S^{-1}(a_{(2, 3)}b_{(3, 3)})
a_{(2, 1)}b_{(3, 1)})\sigma (a_{(2, 2)}b_{(3, 2)}, c_5)\\
&=&A(1)B(S^{-1}(a_5)a_1)C(S^{-1}(a_4b_{11})a_2b_9S^{-1}(b_8)b_6
S^{-1}(b_5)b_1)\\
&&D(S^{-1}(c_3)S^{-1}(b_4)b_2c_1)\sigma (b_3, c_2d)\sigma ^{-1}(b_7, c_4)
\sigma (a_3b_{10}, c_5)\\
&=&A(1)B(S^{-1}(a_5)a_1)C(S^{-1}(b_7)S^{-1}(a_4)a_2b_1)
D(S^{-1}(b_4c_3)b_2c_1)\\
&&\sigma (b_3, c_2d)\sigma ^{-1}(b_5, c_4)\sigma (a_3b_6, c_5)\\
{\rm (\ref{lazy2})}&=&A(1)B(S^{-1}(a_5)a_1)C(S^{-1}(b_7)S^{-1}(a_4)a_2b_1)
D(S^{-1}(b_4c_3)b_3c_2d_3S^{-1}(d_1))\\
&&\sigma (b_2, c_1d_2)\sigma ^{-1}(b_5, c_4)\sigma (a_3b_6, c_5)\\
&=&A(1)B(S^{-1}(a_5)a_1)C(S^{-1}(b_5)S^{-1}(a_4)a_2b_1)
D(d_3S^{-1}(d_1))\\
&&\sigma (b_2, c_1d_2)\sigma ^{-1}(b_3, c_2)\sigma (a_3b_4, c_3),
\end{eqnarray*}
and we see that the two terms are equal.\\
Assume now that $\sigma \in Z^2_{NL}(H)$; we prove (\ref{verypure1}) for 
$\overline{\sigma }$ and the elements $A\ot a$, $B\ot b$, $C\ot c$ in 
$D(H)$. We compute:\\[2mm]
${\;\;\;\;\;\;\;\;}$$\overline{\sigma }(A\ot a, (B\ot b)_1)
\overline{\sigma }((B\ot b)_2, C\ot c)$
\begin{eqnarray*}
&=&\overline{\sigma }(A\ot a, B_2\ot b_1)
\overline{\sigma }(B_1\ot b_2, C\ot c)\\
{\rm (\ref{exti})}&=&A(1)B_2(S^{-1}(a_3)a_1)\sigma (a_2, b_1)
B_1(1)C(S^{-1}(b_{(2, 3)})b_{(2, 1)})\sigma (b_{(2, 2)}, c)\\
&=&A(1)B(S^{-1}(a_3)a_1)
C(S^{-1}(b_4)b_2)\sigma (a_2, b_1)\sigma (b_3, c)\\
&=&A(1)B(S^{-1}(a_3)a_1)
C(S^{-1}(b_4)b_2c_2S^{-1}(c_1))\sigma (a_2, b_1)\sigma (b_3, c_3)\\
{\rm (\ref{lazy2})}&=&A(1)B(S^{-1}(a_3)a_1)
C(S^{-1}(b_4)b_3c_3S^{-1}(c_1))\sigma (a_2, b_1)\sigma (b_2, c_2)\\
{\rm (\ref{verypure1})}&=&A(1)B(S^{-1}(a_3)a_1)
C(c_3S^{-1}(c_1))\sigma (b_1, c_2)\sigma (a_2, b_2),
\end{eqnarray*}
${\;\;\;\;\;\;\;\;}$$\overline{\sigma }((B\ot b)_1, C\ot c)
\overline{\sigma }(A\ot a, (B\ot b)_2)$
\begin{eqnarray*}
&=&\overline{\sigma }(B_2\ot b_1, C\ot c)
\overline{\sigma }(A\ot a, B_1\ot b_2)\\
{\rm (\ref{exti})}&=&B_2(1)C(S^{-1}(b_{(1, 3)})b_{(1, 1)})
\sigma (b_{(1, 2)}, c)A(1)B_1(S^{-1}(a_3)a_1)\sigma (a_2, b_2)\\
&=&A(1)B(S^{-1}(a_3)a_1)C(S^{-1}(b_3)b_1)\sigma (b_2, c)\sigma (a_2, b_4)\\
&=&A(1)B(S^{-1}(a_3)a_1)C(S^{-1}(b_3)b_1c_2S^{-1}(c_1))
\sigma (b_2, c_3)\sigma (a_2, b_4)\\
{\rm (\ref{lazy2})}&=&A(1)B(S^{-1}(a_3)a_1)C(S^{-1}(b_3)b_2c_3S^{-1}(c_1))
\sigma (b_1, c_2)\sigma (a_2, b_4)\\
&=&A(1)B(S^{-1}(a_3)a_1)C(c_3S^{-1}(c_1))
\sigma (b_1, c_2)\sigma (a_2, b_2), 
\end{eqnarray*}   
finishing the proof.
\end{proof}
From a similar computation, the following result follows. 
\begin{proposition}
Let $\gamma \in Reg^1_L(H)$ and define
\begin{eqnarray*}
&&\overline{\gamma }:D(H)\rightarrow k, \;\;\overline{\gamma }(p\ot h)=
p(1)\gamma (h), \;\;\forall \;\;p\in H^*, \;h\in H.
\end{eqnarray*}
If $\gamma \in Reg^1_{PL}(H)$ then $\overline{\gamma }\in Reg^1_{PL}(D(H))$ 
and if $\gamma \in Reg^1_{NL}(H)$ then $\overline{\gamma }\in 
Reg^1_{NL}(D(H))$.
\end{proposition}
We recall now from \cite{rad} the construction of the  
Radford biproduct. Let $H$ be a bialgebra and $B$ a vector space
such that $(B, 1_B)$ is an algebra (with multiplication denoted by
$b\ot c\mapsto bc$ for all $b, c\in B$) and $(B, \Delta _B,
\varepsilon _B)$ is a coalgebra. The pair $(H, B)$ is called {\it
admissible} if $B$ is endowed with a left $H$-module structure
(denoted by $h\ot b\mapsto h\cdot b)$ and with a left
$H$-comodule structure (denoted by $b\mapsto b^{(-1)}\ot
b^{(0)}\in H\ot B)$ such that:\\
(1) $B$ is a left $H$-module algebra;

(2) $B$ is a left $H$-comodule algebra;

(3) $B$ is a left $H$-comodule coalgebra, that is, for all $b\in
B$:
\begin{eqnarray}
&&b_1^{(-1)}b_2^{(-1)}\ot b_1^{(0)}\ot b_2^{(0)}=b^{(-1)}\ot (b^{(0)})_1
\ot (b^{(0)})_2, \label{r1} \\
&&b^{(-1)}\varepsilon _B(b^{(0)})=\varepsilon _B(b)1_H. \label{r2}
\end{eqnarray}

(4) $B$ is a left $H$-module coalgebra, that is, for all $h\in H$
and $b\in B$:
\begin{eqnarray}
&&\Delta _B(h\cdot b)=h_1\cdot b_1\ot h_2\cdot b_2, \label{r3} \\
&&\varepsilon _B(h\cdot b)=\varepsilon _H(h)\varepsilon _B(b). \label{r4}
\end{eqnarray}

(5) $\varepsilon _B$ is an algebra map and $\Delta _B(1_B)=1_B\ot
1_B$;

(6) The following relations hold for all $h\in H$ and $b, c\in B$:
\begin{eqnarray}
&&\Delta _B(bc)=b_1(b_2^{(-1)}\cdot c_1)\ot b_2^{(0)}c_2, \label{r5} \\
&&(h_1\cdot b)^{(-1)}h_2\ot (h_1\cdot b)^{(0)}=h_1b^{(-1)}\ot
h_2\cdot b^{(0)}. \label{r6}
\end{eqnarray}
If $(H, B)$ is an admissible pair, then we know from \cite{rad}
that the smash product algebra structure and smash coproduct
coalgebra structure on $B\ot H$ afford $B\ot H$ a bialgebra
structure, denoted by $B\times H$ and called the {\it smash
biproduct} or {\it Radford biproduct}. Its comultiplication is
given by
\begin{eqnarray}
&&\Delta (b\times h)=(b_1\times b_2^{(-1)}h_1)\ot (b_2^{(0)}\times h_2),
\label{r7}
\end{eqnarray}
for all $b\in B$, $h\in H$, and its counit is $\varepsilon _B\ot
\varepsilon _H$. If $H$ is a Hopf algebra with antipode 
$S_H$ and $(H, B)$ is an admissible pair such that there exists
$S_B\in Hom (B, B)$ a convolution inverse for $id_B$, then
$B\times H$ is a Hopf algebra with antipode
\begin{eqnarray}
S(b\times h)=(1\times S_H(b^{(-1)}h))(S_B(b^{(0)})\times 1), \label{r8}
\end{eqnarray}
for all $h\in H$, $b\in B$. In this case, we will say that $(H,
B)$ is a {\it Hopf admissible pair}. For a Hopf algebra $H$, it is
well-known (see \cite{m}, \cite{mon}) that   
$(H, B)$ being an admissible pair (respectively Hopf admissible
pair) is equivalent to $B$ being a bialgebra (respectively Hopf
algebra) in the Yetter-Drinfeld category $_H^H{\cal
YD}$.\\[2mm]
Let ${\cal C}$ be a braided monoidal category and $B$ a Hopf
algebra in ${\cal C}$. Then, just as if $B$ would be a usual Hopf
algebra, one can define 2-cocycles, crossed products, Galois
extensions, etc, for $B$ in ${\cal C}$, see for instance
\cite{bespalov}, \cite{zhang}. Also, one can define lazy 
2-cocycles, lazy 2-coboundaries and the second lazy cohomology
group $H^2_L(B)=Z^2_L(B)/B^2_L(B)$, see \cite{cp}. We recall these concepts   
in the case when ${\cal C}$=$_H^H{\cal YD}$, the
category of left Yetter-Drinfeld modules over a Hopf algebra $H$,
and $B$ a Hopf algebra in $_H^H{\cal YD}$ (that is, $(H, B)$ is a
Hopf admissible pair, so $B\times H$ is a Hopf algebra).\\
If $M, N\in {}_H^H{\cal YD}$, then $M\ot N\in {}_H^H{\cal YD}$
with module structure $h\cdot (m\ot n)= h_1\cdot m\ot h_2\cdot n$
and comodule structure $m\ot n\mapsto m_{<-1>}n_{<-1>}\ot
(m_{<0>}\ot n_{<0>})$, where $m\mapsto m_{<-1>}\ot m_{<0>}$ and
$n\mapsto n_{<-1>}\ot n_{<0>}$ are the comodule structures of $M$
and $N$, and the braiding is given by
\begin{eqnarray}
&&c_{M, N}:M\ot N\rightarrow N\ot M, \;\;\;\;c_{M, N}(m\ot n)=m_{<-1>}
\cdot n\ot m_{<0>}. \label{b1}
\end{eqnarray}
Hence, the coalgebra structure of $B\ot B$ in $_H^H{\cal YD}$ is given by
\begin{eqnarray*}
\Delta _{B\ot B}(b\ot b')&=&(id \ot c_{B, B}\ot id)\circ (\Delta _B\ot
\Delta _B)(b\ot b')\\
&=&(b_1\ot b_2^{(-1)}\cdot b'_1)\ot (b_2^{(0)}\ot b'_2).
\end{eqnarray*}
So, if $\sigma ,\tau:B\ot B\rightarrow k$ are morphisms in
$_H^H{\cal YD}$, their convolution in $_H^H{\cal YD}$ is given by:
\begin{eqnarray}
&&(\sigma * \tau)(b\ot b')=\sigma (b_1\ot b_2^{(-1)}\cdot b'_1)
\tau (b_2^{(0)}\ot b'_2). \label{b3}
\end{eqnarray}
Let $\sigma :B\ot B\rightarrow k$ be a morphism in
$_H^H{\cal YD}$, that is, it satisfies the conditions:
\begin{eqnarray}
&&\sigma (h_1\cdot b\ot h_2\cdot b')=\varepsilon (h)\sigma (b\ot b'),
\label{b4} \\
&&\sigma (b^{(0)}\ot b^{'(0)})b^{(-1)}b^{'(-1)}=\sigma (b\ot b')1_H,
\label{b5}
\end{eqnarray}
for all $h\in H$ and $b, b'\in B$. Then $\sigma $ is a lazy element if it
satisfies the categorical laziness condition (for all $b, b'\in B$):
\begin{eqnarray}
&&\sigma (b_1\ot b_2^{(-1)}\cdot b'_1)b_2^{(0)}b'_2=\sigma (b_2^{(0)}\ot
b'_2)b_1(b_2^{(-1)}\cdot b'_1). \label{b6}
\end{eqnarray}
Let $\sigma :B\ot B\rightarrow k$ be a normalized left 2-cocycle in
$_H^H{\cal YD}$, that is $\sigma $ is a normalized morphism in
$_H^H{\cal YD}$ satisfying the categorical left 2-cocycle condition
\begin{eqnarray}
\sigma (a_1\ot a_2^{(-1)}\cdot b_1)\sigma (a_2^{(0)}b_2\ot c)=
\sigma (b_1\ot b_2^{(-1)}\cdot c_1)\sigma (a\ot b_2^{(0)}c_2), \label{b6.5}
\end{eqnarray}
for all $a, b, c\in B$. Then we can consider the crossed product
$_{\sigma }B=k\# _{\sigma }B$ as in \cite{zhang}, which is an algebra in
$_H^H{\cal YD}$, and whose multiplication is:
\begin{eqnarray}
&&b\cdot b'=\sigma (b_1\ot b_2^{(-1)}\cdot b'_1)b_2^{(0)}b'_2.\label{b7}
\end{eqnarray}
Since $_{\sigma }B$ is an algebra in $_H^H{\cal YD}$, it is in
particular a left $H$-module algebra, so one can consider the 
smash product $_{\sigma }B\# H$. \\
Let $\gamma :B\rightarrow k$ be a morphism in $_H^H{\cal 
YD}$, that is
\begin{eqnarray}
&&\gamma (h\cdot b)=\varepsilon (h)\gamma (b), \label{b8}\\
&&\gamma (b^{(0)})b^{(-1)}=\gamma (b)1_H, \label{b9}
\end{eqnarray}
for all $h\in H$ and $b\in B$. If $\gamma $ is normalized and
convolution invertible in $_H^H{\cal YD}$, with convolution
inverse $\gamma ^{-1}$ in $_H^H{\cal YD}$, the analogue of the
operator $D^1$ is given in $_H^H{\cal YD}$ by:
\begin{eqnarray*}
D^1(\gamma )(b\ot b')&=&\gamma (b_1)\gamma (b_2^{(-1)}\cdot b'_1)
\gamma ^{-1}(b_2^{(0)}b'_2)  \\
{\rm (\ref{b8})}&=&\gamma (b_1)\gamma (b'_1)\gamma ^{-1}(b_2b'_2),  
\end{eqnarray*}
that is $D^1$ is given by the same formula as for ordinary Hopf
algebras. For a morphism $\gamma :B\rightarrow k$ in $_H^H{\cal
YD}$, the laziness condition is identical to the usual one:
$\gamma (b_1)b_2=b_1\gamma (b_2)$ for all $b\in B$.\\
We recall also the following result from \cite{cp}.
\begin{theorem}(\cite{cp})\label{main}
Let $(H,B)$ be a Hopf admissible pair. \vspace{-5pt}
\begin{enumerate}
\itemsep 0pt
\item[(i)] For a normalized left 2-cocycle $\sigma :B\ot B\rightarrow
k$ in $_H^H{\cal YD}$ define $\overline{\sigma }:(B\times H)\ot
(B\times H) \rightarrow k$,
\begin{eqnarray}
&&\overline{\sigma }(b\times h, b'\times h')=\sigma (b\ot h\cdot b')
\varepsilon (h'). \label{v1}
\end{eqnarray}
Then $\overline{\sigma }$ is a normalized left 2-cocycle on
$B\times H$ and we have $_{\sigma }B\# H=$$\;_{\overline{\sigma
}}(B\times H)$ as algebras. Moreover, $\overline{\sigma }$ is
unique with this property.

\item[(ii)] If $\sigma $ is convolution invertible in $_H^H{\cal
YD}$, then $\overline{\sigma }$ is convolution invertible, with
inverse
\begin{eqnarray}
&&\overline{\sigma }^{-1}(b\times h, b'\times h')=
\sigma ^{-1}(b\ot h\cdot b')
\varepsilon (h'), \label{v2}
\end{eqnarray}
where $\sigma ^{-1}$ is the convolution inverse of $\sigma $ in
$_H^H{\cal YD}$.

\item[(iii)] If $\sigma $ is lazy in $_H^H{\cal YD}$, then
$\overline{\sigma }$ is lazy.

\item[(iv)] If $\sigma,\tau :B\ot B\rightarrow k$ are lazy 2-cocycles
in $_H^H{\cal YD}$, then $\overline{\sigma *\tau
}=\overline{\sigma } *\overline{\tau}$, hence the map $\sigma
\mapsto \overline{\sigma }$ is a group homomorphism from
$Z^2_L(B)$ to $Z^2_L(B\times H)$.

\item[(v)] If $\gamma :B\rightarrow k$ is a normalized and
convolution invertible morphism in $_H^H{\cal YD}$, define
$\overline{\gamma }:B\times H \rightarrow k$ by
\begin{eqnarray}
&&\overline{\gamma }(b\times h)=\gamma (b)\varepsilon (h). \label{v3}
\end{eqnarray}
Then $\overline{\gamma }$ is normalized and convolution
invertible and $\overline{D^1(\gamma )}=D^1(\overline{\gamma })$.
If $\gamma $ is lazy in $_H^H{\cal YD}$, then $\overline{\gamma
}$ is also lazy.

\item[(vi)] If $\sigma $ is a lazy 2-coboundary for $B$ in $_H^H{\cal
YD}$, then $\overline{\sigma }$ is a lazy 2-coboundary for
$B\times H$, so the group homomorphism $Z^2_L(B)\rightarrow
Z^2_L(B\times H)$, $\sigma \mapsto \overline{\sigma }$,
factorizes to a group homomorphism $H^2_L(B)\rightarrow
H^2_L(B\times H)$.
\end{enumerate}
\end{theorem}
For a morphism $\sigma :B\ot B\rightarrow k$ in $_H^H{\cal YD}$, we 
record the following useful formula  
\begin{eqnarray}
&&\sigma (a\ot h\cdot b)=\sigma (S^{-1}(h)\cdot a\ot b), \label{consmor}
\end{eqnarray}
for all $a, b\in B$ and $h\in H$, which is obtained as follows:
\begin{eqnarray*}
\sigma (a\ot h\cdot b)&=&\sigma (h_2S^{-1}(h_1)\cdot a \ot h_3\cdot b)\\
{\rm (\ref{b4})}&=&\sigma (S^{-1}(h)\cdot a\ot b).
\end{eqnarray*}
As for the 2-cocycle condition and laziness condition, there exists a 
categorical analogue for the purity condition (\ref{pure1}), 
which is obtained by appropriately introducing the braiding in 
(\ref{pure1}); for ${\cal C}=$$\;_H^H{\cal YD}$, the condition which is   
obtained may be simplified using repeatedly the formulae  
(\ref{r1}), (\ref{b4}), (\ref{r6}), (\ref{b5}), so we arrive at the 
following concept: 
\begin{definition}{\rm 
Let $(H, B)$ be a Hopf admissible pair and $\sigma :B\ot B\rightarrow k$ 
a normalized and convolution invertible morphism in 
$_H^H{\cal YD}$, with convolution inverse $\sigma ^{-1}$ in 
$_H^H{\cal YD}$. We call $\sigma $ {\it pure} in $_H^H{\cal YD}$ if it 
satisfies the condition (for all $a, b, c, d\in B$):
\begin{eqnarray*}
&&\sigma (ab_1\ot b_2^{(-1)}\cdot c_1)\sigma ^{-1}((b_2^{(0)})_1\ot   
(b_2^{(0)})_2^{(-1)}\cdot c_2)\sigma ((b_2^{(0)})_2^{(0)}\ot c_3d)
\end{eqnarray*}
\begin{eqnarray*}
&&\;\;\;\;\;\;\;\;\;\;\;\;\;\;\;\;\;
=\sigma (b_1\ot [b_2^{(-1)}\cdot c_1][(b_2^{(0)})_2^{(0)(-1)}c_3^{(-1)}
\cdot d])\sigma ^{-1}((b_2^{(0)})_1\ot (b_2^{(0)})_2^{(-1)}\cdot c_2)
\end{eqnarray*}
\begin{eqnarray}
&&\sigma (a(b_2^{(0)})_2^{(0)(0)}\ot c_3^{(0)}). \label{pure4}
\end{eqnarray}
}
\end{definition} 
Similarly, we have the categorical analogue of the condition 
(\ref{verypure1}). 
\begin{definition}{\rm 
Let $(H, B)$ be a Hopf admissible pair and $\sigma :B\ot B\rightarrow k$ 
a normalized and convolution invertible morphism in 
$_H^H{\cal YD}$. We call $\sigma $ {\it neat} in $_H^H{\cal YD}$ if it    
satisfies the condition:
\begin{eqnarray}
\sigma (a\ot b_1)\sigma (b_2\ot c)=\sigma (a^{(-1)}\cdot b_1\ot c)
\sigma (a^{(0)}\ot b_2), \;\;\; \forall \;a, b, c\in B. 
 \label{verypure2}
\end{eqnarray}
}
\end{definition}
\begin{remark}{\em 
It is not straightforward to prove that a neat lazy cocycle   
$\sigma :B\ot B\rightarrow k$ in $_H^H{\cal YD}$ is pure in $_H^H{\cal YD}$. 
We will see an indirect proof below.}
\end{remark}
Motivated by Theorem \ref{main}, we prove the following result. 
\begin{theorem}\label{pvp}
Let $(H, B)$ be a Hopf admissible pair and $\sigma :B\ot B\rightarrow k$ 
pure (respectively neat) in $_H^H{\cal YD}$. If we define 
$\overline{\sigma }:(B\times H)\ot (B\times H)\rightarrow k$ by formula 
(\ref{v1}), then $\overline{\sigma }$ is pure (respectively neat).  
In particular, if $\sigma $ is a pure (respectively neat) lazy cocycle  
in $_H^H{\cal YD}$, then $\overline{\sigma }$ is a pure (respectively 
neat) lazy cocycle for $B\times H$. 
\end{theorem}
\begin{proof}
Note first that $\overline{\sigma }$ is convolution invertible, with 
convolution inverse given by (\ref{v2}) ($\sigma $ does not have to be a 
2-cocycle for this). Now let $a, b, c, d\in B$ and 
$h, g, l, t\in H$ and assume that $\sigma $ is pure in $_H^H{\cal YD}$; 
we prove the purity condition (\ref{pure1}) for $\overline{\sigma }$ 
on $B\times H$, for the elements $a\times h$, $b\times g$, $c\times l$,  
$d\times t$. First we compute the right hand side of (\ref{pure1}):\\[2mm]
${\;\;\;}$
$\overline{\sigma }((b\times g)_1, (c\times l)_1(d\times t))
\overline{\sigma }^{-1}((b\times g)_2, (c\times l)_2)\overline{\sigma }
((a\times h)(b\times g)_3, (c\times l)_3)$
\begin{eqnarray*}
{\rm (\ref{r7})}&=&\overline{\sigma }(b_1\times b_2^{(-1)}g_1, 
c_1((c_2^{(-1)})_1l_1\cdot d)\times (c_2^{(-1)})_2l_2t)\\
&&\overline{\sigma }^{-1}((b_2^{(0)})_1\times (b_2^{(0)})_2^{(-1)}g_2, 
(c_2^{(0)})_1\times (c_2^{(0)})_2^{(-1)}l_3)\\
&&\overline{\sigma }(a(h_1\cdot (b_2^{(0)})_2^{(0)})\times h_2g_3, 
(c_2^{(0)})_2^{(0)}\times l_4)\\
{\rm (\ref{v1}, \ref{v2})}&=&\sigma (b_1\ot [(b_2^{(-1)})_1g_1\cdot c_1]
[(b_2^{(-1)})_2g_2c_2^{(-1)}l\cdot d])\\
&&\sigma ^{-1}((b_2^{(0)})_1\ot (b_2^{(0)})_2^{(-1)}g_3\cdot 
(c_2^{(0)})_1)\\
&&\sigma (a(h_1\cdot (b_2^{(0)})_2^{(0)})\ot h_2g_4\cdot (c_2^{(0)})_2)
\varepsilon (t)\\
{\rm (\ref{r1})}&=&\sigma (b_1\ot [(b_2^{(-1)})_1g_1\cdot c_1]
[(b_2^{(-1)})_2g_2c_2^{(-1)}c_3^{(-1)}l\cdot d])\\
&&\sigma ^{-1}((b_2^{(0)})_1\ot (b_2^{(0)})_2^{(-1)}g_3\cdot 
c_2^{(0)})\\
&&\sigma (a(h_1\cdot (b_2^{(0)})_2^{(0)})\ot h_2g_4\cdot c_3^{(0)})
\varepsilon (t)\\
{\rm (\ref{r1})}&=&\sigma (b_1\ot [(b_2^{(-1)})_1(b_3^{(-1)})_1
g_1\cdot c_1]
[(b_2^{(-1)})_2(b_3^{(-1)})_2
g_2c_2^{(-1)}c_3^{(-1)}l\cdot d])\\
&&\sigma ^{-1}(b_2^{(0)}\ot (b_3^{(-1)})_3g_3\cdot 
c_2^{(0)})\\
&&\sigma (a(h_1\cdot b_3^{(0)})\ot h_2g_4\cdot c_3^{(0)})
\varepsilon (t)\\
{\rm (\ref{r6})}&=&\sigma (b_1\ot [b_2^{(-1)}(b_3^{(-1)})_1
g_1\cdot c_1]
[b_2^{(0)(-1)}((b_3^{(-1)})_2
g_2\cdot c_2)^{(-1)}(b_3^{(-1)})_3g_3
c_3^{(-1)}l\cdot d])\\
&&\sigma ^{-1}(b_2^{(0)(0)}\ot ((b_3^{(-1)})_2g_2\cdot  
c_2)^{(0)})\\
&&\sigma (a(h_1\cdot b_3^{(0)})\ot h_2g_4\cdot c_3^{(0)})
\varepsilon (t)\\
{\rm (\ref{b5})}&=&\sigma (b_1\ot [b_2^{(-1)}b_3^{(-1)}
g_1\cdot c_1]
[(b_3^{(0)(-1)})_2g_3
c_3^{(-1)}l\cdot d])\\
&&\sigma ^{-1}(b_2^{(0)}\ot (b_3^{(0)(-1)})_1g_2\cdot   
c_2)\\
&&\sigma (a(h_1\cdot b_3^{(0)(0)})\ot h_2g_4\cdot c_3^{(0)})
\varepsilon (t)\\
{\rm (\ref{r1})}&=&\sigma (b_1\ot [b_2^{(-1)}
g_1\cdot c_1]
[((b_2^{(0)})_2^{(-1)})_2
g_3c_3^{(-1)}l\cdot d])\\
&&\sigma ^{-1}((b_2^{(0)})_1\ot ((b_2^{(0)})_2^{(-1)})_1g_2\cdot    
c_2)\\
&&\sigma (a(h_1\cdot (b_2^{(0)})_2^{(0)})\ot h_2g_4\cdot c_3^{(0)})
\varepsilon (t).
\end{eqnarray*}
Now we compute the left hand side of (\ref{pure1}):\\[2mm]
${\;\;\;}$$\overline{\sigma }((a\times h)(b\times g)_1, (c\times l)_1)
\overline{\sigma }^{-1}((b\times g)_2, (c\times l)_2)
\sigma ((b\times g)_3, (c\times l)_3(d\times t))$
\begin{eqnarray*}
{\rm (\ref{r7})}&=&\overline{\sigma }(a(h_1\cdot b_1)\times 
h_2b_2^{(-1)}g_1, c_1\times c_2^{(-1)}l_1)\\
&&\overline{\sigma }^{-1}((b_2^{(0)})_1\times (b_2^{(0)})_2^{(-1)}g_2, 
(c_2^{(0)})_1\times (c_2^{(0)})_2^{(-1)}l_2)\\
&&\overline{\sigma }((b_2^{(0)})_2^{(0)}\times g_3, (c_2^{(0)})_2^{(0)}
(l_3\cdot d)\times l_4t)\\
{\rm (\ref{v1}, \ref{v2})}&=&\sigma (a(h_1\cdot b_1)\ot h_2b_2^{(-1)}
g_1\cdot c_1)\\
&&\sigma ^{-1}((b_2^{(0)})_1\ot (b_2^{(0)})_2^{(-1)}g_2\cdot c_2)\\
&&\sigma ((b_2^{(0)})_2^{(0)}\ot [g_3\cdot c_3][g_4l\cdot d])
\varepsilon (t)\\
&=&\sigma (a(h_1\cdot b_1)\ot h_2b_2^{(-1)}
g_1\cdot c_1)\\
&&\sigma ^{-1}((S^{-1}(h_4)h_3\cdot b_2^{(0)})_1
\ot (S^{-1}(h_4)h_3\cdot b_2^{(0)})_2^{(-1)}g_2\cdot c_2)\\
&&\sigma ((S^{-1}(h_4)h_3\cdot b_2^{(0)})_2^{(0)}
\ot [g_3\cdot c_3][g_4l\cdot d])\varepsilon (t)\\
{\rm (\ref{r3})}&=&\sigma (a(h_1\cdot b_1)\ot h_2b_2^{(-1)}
g_1\cdot c_1)\\
&&\sigma ^{-1}(S^{-1}(h_4)_1\cdot (h_3\cdot b_2^{(0)})_1
\ot [S^{-1}(h_4)_2\cdot (h_3\cdot b_2^{(0)})_2]^{(-1)}g_2\cdot c_2)\\
&&\sigma ([S^{-1}(h_4)_2\cdot (h_3\cdot b_2^{(0)})_2]^{(0)}
\ot [g_3\cdot c_3][g_4l\cdot d])\varepsilon (t)\\
{\rm (\ref{r6})}&=&\sigma (a(h_1\cdot b_1)\ot (h_2\cdot b_2)^{(-1)}
h_3g_1\cdot c_1)\\
&&\sigma ^{-1}(S^{-1}(h_4)_1\cdot [(h_2\cdot b_2)^{(0)}]_1
\ot [S^{-1}(h_4)_2\cdot [(h_2\cdot b_2)^{(0)}]_2]^{(-1)}g_2\cdot c_2)\\
&&\sigma ([S^{-1}(h_4)_2\cdot [(h_2\cdot b_2)^{(0)}]_2]^{(0)}
\ot [g_3\cdot c_3][g_4l\cdot d])\varepsilon (t)\\
{\rm (\ref{r6})}&=&\sigma (a(h_1\cdot b_1)\ot (h_2\cdot b_2)^{(-1)}
h_3g_1\cdot c_1)\\
&&\sigma ^{-1}(S^{-1}(h_4)_1\cdot ((h_2\cdot b_2)^{(0)})_1
\ot S^{-1}(h_4)_2((h_2\cdot b_2)^{(0)})_2^{(-1)}S(S^{-1}(h_4)_4)
g_2\cdot c_2)\\
&&\sigma (S^{-1}(h_4)_3\cdot ((h_2\cdot b_2)^{(0)})_2^{(0)}
\ot [g_3\cdot c_3][g_4l\cdot d])\varepsilon (t)\\
{\rm (\ref{b4})}&=&\sigma (a(h_1\cdot b_1)\ot (h_2\cdot b_2)^{(-1)}
h_3g_1\cdot c_1)\\
&&\sigma ^{-1}(((h_2\cdot b_2)^{(0)})_1
\ot ((h_2\cdot b_2)^{(0)})_2^{(-1)}h_4g_2\cdot c_2)\\
&&\sigma (S^{-1}(h_5)\cdot ((h_2\cdot b_2)^{(0)})_2^{(0)}
\ot [g_3\cdot c_3][g_4l\cdot d])\varepsilon (t)\\
{\rm (\ref{r3}, \ref{consmor})}&=&
\sigma (a(h_1\cdot b)_1\ot (h_1\cdot b)_2^{(-1)}h_2g_1\cdot c_1)\\
&&\sigma ^{-1}(((h_1\cdot b)_2^{(0)})_1
\ot ((h_1\cdot b)_2^{(0)})_2^{(-1)}h_3g_2\cdot c_2)\\
&&\sigma (((h_1\cdot b)_2^{(0)})_2^{(0)}
\ot [h_4g_3\cdot c_3][h_5g_4l\cdot d])\varepsilon (t)\\
{\rm (\ref{r3})}&=&
\sigma (a(h_1\cdot b)_1\ot (h_1\cdot b)_2^{(-1)}\cdot (h_2g_1\cdot c)_1)\\
&&\sigma ^{-1}(((h_1\cdot b)_2^{(0)})_1
\ot ((h_1\cdot b)_2^{(0)})_2^{(-1)}\cdot (h_2g_1\cdot c)_2)\\
&&\sigma (((h_1\cdot b)_2^{(0)})_2^{(0)}
\ot (h_2g_1\cdot c)_3(h_3g_2l\cdot d))\varepsilon (t)\\
{\rm (\ref{pure4})}&=&
\sigma ((h_1\cdot b)_1\ot [(h_1\cdot b)_2^{(-1)}\cdot (h_2g_1\cdot c)_1]
[((h_1\cdot b)_2^{(0)})_2^{(0)(-1)}(h_2g_1\cdot c)_3^{(-1)}h_3g_2l\cdot d])\\
&&\sigma ^{-1}(((h_1\cdot b)_2^{(0)})_1
\ot ((h_1\cdot b)_2^{(0)})_2^{(-1)}\cdot (h_2g_1\cdot c)_2)\\
&&\sigma (a((h_1\cdot b)_2^{(0)})_2^{(0)(0)}
\ot (h_2g_1\cdot c)_3^{(0)})\varepsilon (t)\\
{\rm (\ref{r3})}&=&
\sigma (h_1\cdot b_1\ot [(h_2\cdot b_2)^{(-1)}h_3g_1\cdot c_1]
[((h_2\cdot b_2)^{(0)})_2^{(0)(-1)}(h_5g_3\cdot c_3)^{(-1)}h_6g_4l\cdot d])\\
&&\sigma ^{-1}(((h_2\cdot b_2)^{(0)})_1
\ot ((h_2\cdot b_2)^{(0)})_2^{(-1)}h_4g_2\cdot c_2)\\
&&\sigma (a((h_2\cdot b_2)^{(0)})_2^{(0)(0)}
\ot (h_5g_3\cdot c_3)^{(0)})\varepsilon (t)\\
{\rm (\ref{r6})}&=&
\sigma (h_1\cdot b_1\ot [h_2b_2^{(-1)}g_1\cdot c_1]
[(h_3\cdot b_2^{(0)})_2^{(0)(-1)}h_5g_3c_3^{(-1)}l\cdot d])\\
&&\sigma ^{-1}((h_3\cdot b_2^{(0)})_1
\ot (h_3\cdot b_2^{(0)})_2^{(-1)}h_4g_2\cdot c_2)\\
&&\sigma (a(h_3\cdot b_2^{(0)})_2^{(0)(0)}
\ot h_6g_4\cdot c_3^{(0)})\varepsilon (t)\\
{\rm (\ref{r3})}&=&
\sigma (h_1\cdot b_1\ot [h_2b_2^{(-1)}g_1\cdot c_1]
[(h_4\cdot (b_2^{(0)})_2)^{(0)(-1)}h_6g_3c_3^{(-1)}l\cdot d])\\
&&\sigma ^{-1}(h_3\cdot (b_2^{(0)})_1
\ot (h_4\cdot (b_2^{(0)})_2)^{(-1)}h_5g_2\cdot c_2)\\
&&\sigma (a(h_4\cdot (b_2^{(0)})_2)^{(0)(0)}
\ot h_7g_4\cdot c_3^{(0)})\varepsilon (t)\\
{\rm (\ref{r6})}&=&
\sigma (h_1\cdot b_1\ot [h_2b_2^{(-1)}g_1\cdot c_1]
[(h_5\cdot (b_2^{(0)})_2^{(0)})^{(-1)}h_6g_3c_3^{(-1)}l\cdot d])\\
&&\sigma ^{-1}(h_3\cdot (b_2^{(0)})_1
\ot h_4(b_2^{(0)})_2^{(-1)}g_2\cdot c_2)\\
&&\sigma (a(h_5\cdot (b_2^{(0)})_2^{(0)})^{(0)}
\ot h_7g_4\cdot c_3^{(0)})\varepsilon (t)\\
{\rm (\ref{b4})}&=&
\sigma (h_1\cdot b_1\ot [h_2b_2^{(-1)}g_1\cdot c_1]
[(h_3\cdot (b_2^{(0)})_2^{(0)})^{(-1)}h_4g_3c_3^{(-1)}l\cdot d])\\
&&\sigma ^{-1}((b_2^{(0)})_1
\ot (b_2^{(0)})_2^{(-1)}g_2\cdot c_2)\\
&&\sigma (a(h_3\cdot (b_2^{(0)})_2^{(0)})^{(0)}
\ot h_5g_4\cdot c_3^{(0)})\varepsilon (t)\\
{\rm (\ref{r6})}&=&
\sigma (h_1\cdot b_1\ot [h_2b_2^{(-1)}g_1\cdot c_1]
[h_3(b_2^{(0)})_2^{(0)(-1)}g_3c_3^{(-1)}l\cdot d])\\
&&\sigma ^{-1}((b_2^{(0)})_1
\ot (b_2^{(0)})_2^{(-1)}g_2\cdot c_2)\\
&&\sigma (a(h_4\cdot (b_2^{(0)})_2^{(0)(0)})
\ot h_5g_4\cdot c_3^{(0)})\varepsilon (t)\\
{\rm (\ref{b4})}&=&
\sigma (b_1\ot [b_2^{(-1)}g_1\cdot c_1]
[(b_2^{(0)})_2^{(0)(-1)}g_3c_3^{(-1)}l\cdot d])\\
&&\sigma ^{-1}((b_2^{(0)})_1
\ot (b_2^{(0)})_2^{(-1)}g_2\cdot c_2)\\
&&\sigma (a(h_1\cdot (b_2^{(0)})_2^{(0)(0)})
\ot h_2g_4\cdot c_3^{(0)})\varepsilon (t)\\
&=&\sigma (b_1\ot [b_2^{(-1)}g_1\cdot c_1]
[((b_2^{(0)})_2^{(-1)})_2g_3c_3^{(-1)}l\cdot d])\\
&&\sigma ^{-1}((b_2^{(0)})_1
\ot ((b_2^{(0)})_2^{(-1)})_1g_2\cdot c_2)\\
&&\sigma (a(h_1\cdot (b_2^{(0)})_2^{(0)})
\ot h_2g_4\cdot c_3^{(0)})\varepsilon (t), 
\end{eqnarray*}
and we see that the two terms are equal. 
Assume now that $\sigma $ is neat in $_H^H{\cal YD}$; we prove  
(\ref{verypure1}) for $\overline{\sigma }$ on $B\times H$, for the 
elements $a\times h$, $b\times g$, $c\times l$. We compute:\\[2mm]
${\;\;\;\;\;}$$\overline{\sigma }(a\times h, (b\times g)_1)
\overline{\sigma }((b\times g)_2, c\times l)$
\begin{eqnarray*}
&=&\overline{\sigma }(a\times h, b_1\times b_2^{(-1)}g_1)
\overline{\sigma }(b_2^{(0)}\times g_2, c\times l)\\
{\rm (\ref{v1})}&=&\sigma (a\ot h\cdot b_1)
\sigma (b_2\ot g\cdot c)\varepsilon (l)\\
&=&\sigma (a\ot h_1\cdot b_1)
\sigma (S^{-1}(h_3)h_2\cdot b_2\ot g\cdot c)\varepsilon (l)\\
{\rm (\ref{r3})}&=&\sigma (a\ot (h_1\cdot b)_1)
\sigma (S^{-1}(h_2)\cdot (h_1\cdot b)_2\ot g\cdot c)\varepsilon (l)\\
{\rm (\ref{consmor})}&=&\sigma (a\ot (h_1\cdot b)_1)
\sigma ((h_1\cdot b)_2\ot h_2g\cdot c)\varepsilon (l)\\
{\rm (\ref{verypure2})}&=&\sigma (a^{(-1)}\cdot (h_1\cdot b)_1\ot 
h_2g\cdot c)\sigma (a^{(0)}\ot (h_1\cdot b)_2)\varepsilon (l)\\
{\rm (\ref{r3})}&=&\sigma (a^{(-1)}h_1\cdot b_1\ot  
h_3g\cdot c)\sigma (a^{(0)}\ot h_2\cdot b_2)\varepsilon (l), 
\end{eqnarray*}
${\;\;\;\;\;}$$\overline{\sigma }((b\times g)_1, c\times l)
\overline{\sigma }(a\times h, (b\times g)_2)$
\begin{eqnarray*}
&=&\overline{\sigma }(b_1\times b_2^{(-1)}g_1, c\times l)
\overline{\sigma }(a\times h, b_2^{(0)}\times g_2)\\
{\rm (\ref{v1})}&=&\sigma (b_1\ot b_2^{(-1)}g\cdot c)\sigma (a\ot 
h\cdot b_2^{(0)})\varepsilon (l)\\
&=&\sigma (b_1\ot S(h_1)h_2b_2^{(-1)}g\cdot c)\sigma (a\ot 
h_3\cdot b_2^{(0)})\varepsilon (l)\\
{\rm (\ref{r6})}&=&\sigma (b_1\ot S(h_1)(h_2\cdot b_2)^{(-1)}
h_3g\cdot c)\sigma (a\ot  
(h_2\cdot b_2)^{(0)})\varepsilon (l)\\
&=&\sigma (b_1\ot S(h_1)S((a^{(-1)})_1)(a^{(-1)})_2
(h_2\cdot b_2)^{(-1)}h_3g\cdot c)\\
&&\sigma (a^{(0)}\ot    
(h_2\cdot b_2)^{(0)})\varepsilon (l)\\
&=&\sigma (b_1\ot S(a^{(-1)}h_1)a^{(0)(-1)}
(h_2\cdot b_2)^{(-1)}h_3g\cdot c)\\
&&\sigma (a^{(0)(0)}\ot     
(h_2\cdot b_2)^{(0)})\varepsilon (l)\\
{\rm (\ref{b5}, \ref{consmor})}&=&
\sigma (a^{(-1)}h_1\cdot b_1\ot 
h_3g\cdot c)\sigma (a^{(0)}\ot     
h_2\cdot b_2)\varepsilon (l),
\end{eqnarray*}
and the proof is finished.
\end{proof}  
\begin{remark}{\em \label{remarca}
Let $(H, B)$ be a Hopf admissible pair and $\sigma :B\ot B\rightarrow k$ 
a normalized and convolution invertible morphism in 
$_H^H{\cal YD}$, and define $\overline{\sigma }$ by formula 
(\ref{v1}). From the computation in the proof of Theorem \ref{pvp} it 
follows that, conversely, if $\overline{\sigma }$ is pure 
(respectively neat) on $B\times H$, then $\sigma $ is pure  
(respectively neat) in $_H^H{\cal YD}$. Together with  
Proposition \ref{very} and Theorems \ref{main} and \ref{pvp} 
this proves that, if   
$\sigma $ is a neat lazy cocycle  
in $_H^H{\cal YD}$, then $\sigma $ is a pure lazy cocycle 
in $_H^H{\cal YD}$.}
\end{remark}
There exist also categorical analogues of the relations (\ref{pure0}) and 
(\ref{pure2}); the one corresponding to (\ref{pure0}) looks very 
complicated, so we treat only the analogue of (\ref{pure2}). 
\begin{definition}{\em 
Let $(H, B)$ be a Hopf admissible pair and $\gamma :B\rightarrow k$ a 
normalized and convolution invertible morphism in $_H^H{\cal YD}$; we 
call $\gamma $ {\it neat} in $_H^H{\cal YD}$ if it satisfies the condition 
\begin{eqnarray}
&&\gamma (ab_1)\gamma (b_2c)=\gamma (b_1(b_2^{(-1)}\cdot c))
\gamma (ab_2^{(0)}),\;\;\;\;\forall \; a, b, c\in B. 
\label{verypure3}
\end{eqnarray}
} 
\end{definition}
\begin{proposition}
If $\gamma $ is neat in $_H^H{\cal YD}$, then the map $\overline{\gamma }:
B\times H\rightarrow k$ given by (\ref{v3}) is neat.
\end{proposition}
\begin{proof}
Let $a, b, c\in B$ and $h, g, l\in H$; we check (\ref{pure2}) for the 
elements $a\times h$, $b\times g$, $c\times l$. We compute:\\[2mm]
${\;\;\;\;\;\;\;\;}$$\overline{\gamma }((a\times h)(b\times g)_1)
\overline{\gamma }((b\times g)_2(c\times l))$
\begin{eqnarray*}
{\rm (\ref{r7})}&=&\overline{\gamma }((a\times h)(b_1\times b_2^{(-1)}g_1))
\overline{\gamma }((b_2^{(0)}\times g_2)(c\times l))\\
&=&\overline{\gamma }(a(h_1\cdot b_1)\times h_2b_2^{(-1)}g_1)
\overline{\gamma }(b_2^{(0)}(g_2\cdot c)\times g_3l)\\
{\rm (\ref{v3})}&=&\gamma (a(h\cdot b_1))\gamma (b_2(g\cdot c))
\varepsilon (l)\\
&=&\gamma (a(h_1\cdot b_1))\gamma ((S^{-1}(h_3)\cdot (h_2\cdot b_2))
(g\cdot c))\varepsilon (l)\\
{\rm (\ref{r3})}&=&\gamma (a(h_1\cdot b)_1)\gamma ((S^{-1}(h_2)\cdot 
(h_1\cdot b)_2)(g\cdot c))\varepsilon (l)\\
{\rm (\ref{b8})}&=&\gamma (a(h_1\cdot b)_1)\gamma ((h_1\cdot b)_2
(h_2g\cdot c))\varepsilon (l)\\
{\rm (\ref{verypure3})}&=&\gamma ((h_1\cdot b)_1
((h_1\cdot b)_2^{(-1)}h_2g\cdot c))\gamma (a(h_1\cdot b)_2^{(0)})
\varepsilon (l)\\
{\rm (\ref{r3})}&=&\gamma ((h_1\cdot b_1)((h_2\cdot b_2)^{(-1)}h_3g\cdot c))
\gamma (a(h_2\cdot b_2)^{(0)})\varepsilon (l)\\
{\rm (\ref{r6})}&=&\gamma ((h_1\cdot b_1)(h_2b_2^{(-1)}g\cdot c))
\gamma (a(h_3\cdot b_2^{(0)}))\varepsilon (l)\\
{\rm (\ref{b8})}&=&\gamma (b_1(b_2^{(-1)}g\cdot c))\gamma (a(h\cdot  
b_2^{(0)}))\varepsilon (l), 
\end{eqnarray*}
${\;\;\;\;\;\;\;\;}$$\overline{\gamma }((b\times g)_1(c\times l))
\overline{\gamma }((a\times h)(b\times g)_2)$
\begin{eqnarray*}
{\rm (\ref{r7})}&=&\overline{\gamma }((b_1\times b_2^{(-1)}g_1)(c\times l))
\overline{\gamma }((a\times h)(b_2^{(0)}\times g_2))\\
&=&\overline{\gamma }(b_1((b_2^{(-1)})_1g_1\cdot c)\times (b_2^{(-1)})_2g_2l)
\overline{\gamma }(a(h_1\cdot b_2^{(0)})\times h_2g_3)\\
{\rm (\ref{v3})}&=&\gamma (b_1(b_2^{(-1)}g\cdot c))
\gamma (a(h\cdot b_2^{(0)}))\varepsilon (l),
\end{eqnarray*}
and we see that the two terms are equal.
\end{proof}
\begin{remark}{\em 
Combining  Proposition \ref{vpln}, Theorem \ref{main} (v) and Remark 
\ref{remarca}, we obtain: if $\gamma :B\rightarrow k$ is 
neat lazy in $_H^H{\cal YD}$ then $D^1(\gamma )$ is neat lazy in 
$_H^H{\cal YD}$.}  
\end{remark}

\end{document}